\documentclass[12pt]{amsart}
\usepackage{amscd,amsthm,amsmath} 
\newtheorem{thm}{Theorem}[section]
\newtheorem{cor}[thm]{Corollary}
\newtheorem{lem}[thm]{Lemma}
\newtheorem{prop}[thm]{Proposition}

\newtheorem{ques}[thm]{Question}

\theoremstyle{definition}
\newtheorem{defn}[thm]{Definition}
\theoremstyle{remark}
\newtheorem{rem}[thm]{Remark}
\numberwithin{equation}{section}

\begin{document}

\title{Rational homotopy theory and nonnegative curvature}%
\author{Jianzhong Pan }%
\address{Institute of Math.,Academia Sinica ,Beijing 100080, China }%
\email{pjz@math03.math.ac.cn}%

\thanks{The first author is partially supported by the NSFC
projects  10071087 ,
19701032 and ZD9603 of Chinese Academy of Science }%

\date{Aug. 30, 2001}%
\begin{abstract}
  In this note ,  we answer positively a question by Belegradek and
  Kapovitch\cite{BK1} about the relation between rational homotopy theory and
   a problem in Riemannian geometry which asks that total spaces of
   which vector bundles over compact nonnegative curved manifolds admit (complete)
   metrics with nonnegative curvature.
\end{abstract}
\maketitle
 \section{Introduction} \label{S:intro}

Given a Riemannian manifold $M$ with metric $$<  > :TM \times TM
\to TM$$   an affine connection is a bilinear map $$\nabla :
Vec(M) \times Vec(M) \to Vec(M)$$ which satisfies the following
\begin{itemize}
\item{$\nabla_{fV}W = f \nabla_{V}W$}
\item{$\nabla_{V}(fW) = (Vf)W + f\nabla_{V}W$}
\end{itemize}
where $f \in C^{\infty}(M) , V,W \in Vec(M)$

An affine connection is called Levi-Civita connection if it
satisfies also the following
\begin{itemize}
\item{$X<V,W> = <\nabla_{X}V,W>  + <V, \nabla_{X}W>$}
\item{$\nabla_{V}(W) - \nabla_{W}V - [V,W] =0$}
\end{itemize}
where $[V,W]f = (XY - YX)f$ is the Lie bracket.

A fundamental result in Riemannian geometry asserts that
\begin{thm}
For each Riemannian metric, there exists a unique Levi-Civita
connection.
\end{thm}
Given a Riemannian manifold $M$ with Levi-Civita connection ,
there is defined a curvature operator $$R: Vec(M) \times Vec(M)
\times Vec(M) \to  Vec(M)$$ defined by $$R(X,Y)Z =\nabla_X
\nabla_Y Z - \nabla_Y \nabla_X Z -\nabla_{[X,Y] Z} $$ From it one
arrives at an important geometric invariant which is called {\it
Sectional curvature} defined by $$K(\sigma) =
\frac{<R(v,w)w,v>}{<v \wedge w , v \wedge w>}$$ where $\sigma
\subset T_pM$ is a tangent plane at $p \in M$ and $v, w \in
\sigma$ span it. It is well known that $K(\sigma)$ does not depend
on the choice of spanning vectors.

A well known question in Riemannian geometry is
\begin{ques}
Does the restriction on curvature imply the restriction on
topology and vice versa?
\end{ques}
In particular, how does the positive(nonnegative) curvature
restrict the topology of the underlining manifold?

A Riemannian manifold is called positively (or nonnegatively)
curved if,  for any $\sigma$,  $K(\sigma) > 0$   (or $K(\sigma)
\geq 0)$.

For compact manifold, we have the following classical
\begin{thm}
Let $M$ be a compact Riemannian manifold with positive curvature.
Then
\[ \pi_1(M)= \left\{ \begin{array}{ll}
\mbox{finite group} & \mbox{ if $\dim M$ is odd }\\ 0  & \mbox{ if
$\dim M$ is even and  $M$ is orientable} \\ Z_2  & \mbox{ if $\dim
M$ is even and  $M$ is nonorientable}
\end{array} \right.
\]
\end{thm}

The main concern of this note is on noncompact manifold. In this
case there is the following

\begin{thm}
Let $M$ be a complete noncompact Riemannian manifold with
nonnegative curvature. Then  $M$ is diffeomorphic to the total
space of the normal bundle of a compact totally geodesic
submanifold which is called the {\it soul}.
\end{thm}

Another central question in Riemannian geometry is to what extent
the converse is true, or in other words
\begin{ques}
Total spaces of which vector bundles over compact nonnegatively
curved manifolds admit (complete)
   metrics with nonnegative curvature?
\end{ques}

Previously, obstructions to the existence of nonnegatively curved
metrics on vector bundles were only known for a flat soul
~\cite{OW}. No obstructions are known when the soul is
simply-connected. In \cite{BK1} an approach to the reduction of
the problem to the vector bundle over simply connected manifold
was initiated.  The start point is  another result of Cheeger and
Gromoll \cite{CG}  that a finite cover of any closed nonnegatively
curved manifold (throughout the paper by a nonnegatively curved
manifold we mean a complete Riemannian manifold of nonnegative
{\it sectional} curvature) is diffeomorphic to a product of a
torus and a simply-connected closed nonnegatively curved manifold.
It turns out that a similar statement holds for open complete
nonnegatively curved manifolds which is the basis of their
analysis.

\begin{lem}\cite{BK1}\label{mainlemma}
Let $(N,g)$ be a complete nonnegatively curved manifold. Then
there exists a finite cover $N'$ of $N$ diffeomorphic to a product
$M\times T^k$ where $M$ is a complete open simply connected
nonnegatively curved manifold. Moreover, if $S^\prime$ is a soul
of $N'$, then this diffeomorphism can be chosen in such a way that
it takes $S'$ onto $C\times T^k$ where $C$ is a soul of $M$.
\end{lem}
By using this and characteristic classes technique, they proved
that, in  various case, the total spaces of rank $k$ vector
bundles over $C\times T$ admit
 no nonnegatively curved metric if
they do not become the pullback of a bundle over $C$ in a finite
cover. The following is  such an example
\begin{cor}\cite{BK2}
Let $B$ be a closed nonnegatively curved manifold. If $\pi_1(B)$
contains a free abelian subgroup of rank four (two, respectively),
then for each $k\ge 2$ (for $k=2$, respectively) there exists a
finite cover of $B$ over which there exist infinitely many rank
$k$ vector bundles  whose total spaces admit no nonnegatively
curved metrics.
\end{cor}

 Belegradek and  Kapovitch \cite{BK1} are thus lead to the
following
\begin{defn}
Given a closed smooth simply connected manifold $C$, a torus $T$,
and a positive integer $k$, we say that a triple $(C,T,k)$ is {\it
splitting rigid} if any rank $k$ vector bundle over $C\times T$
with nonnegatively curved total space splits, after passing to a
finite cover, as the product of a rank $k$ bundle over $C$ and a
rank zero bundle over $T$.
\end{defn}

Let $\mathcal H$ be the class of simply-connected CW-complexes
whose rational cohomology algebra is finite dimensional, as a
rational vector space, and has no nonzero derivations of negative
degree(see \cite{BK1} for the reason to choose such a class
$\mathcal H$). For example, $\mathcal H$ contains any compact
simply-connected K\"ahler manifold~\cite{Mei2}.

A natural question is
\begin{ques}\cite{BK1}
 Let $C\in \mathcal H$ be a closed smooth manifold.
Is $(C,T,k)$ splitting rigid for any $T$ and $k$?
\end{ques}

The main result in this note is a positive answer to this question
\begin{thm}\label{T:main}
 Let $C\in \mathcal H$ be a closed smooth manifold.
Then $(C,T,k)$ is splitting rigid for any $T$ and $k$.
\end{thm}

In this  paper, all (co)homology groups have rational
coefficients, all manifolds and vector bundles are smooth; all
topological spaces are homotopy equivalent to connected
CW-complexes. $[X,Y]$ will be the based homotopy classes of based
maps between them. $map(X,Y)$ is  the space of maps from $X$ to
$Y$ and $map(X,Y)_f$ is the connected component of $map(X,Y)$
which contains the map $f:X \to Y$.

\section{A splitting criterion}\label{S:phan}

Given a finite cell complex $C$, define $Char(k,C)$ to be the
subspace of $H^*(C)$ which is the direct sum of
$\oplus_{i=1}^{[(k-1)/2]}H^{4i}(C)$ and the subspace equal to
$H^k(C)$ if $k$ is even, and to $H^{4[k/2]}(C)$ if $k$ is odd.
Note that any rational characteristic class of a rank $k$ vector
bundle over $C$ lies in the subalgebra of $H^*(C)$ generated by
$Char(k,C)$.

 Belegradek and  Kapovitch transform the problem of a triple
 $(C,T,k)$ being splitting rigid into a homotopy problem as follows

 \begin{prop}\cite{BK1} \label{characterization}
Let $C$ be a closed simply-connected manifold, $T$ be a torus, $k$
be a positive integer. If any self-homotopy equivalence of
$C\times T$ maps $Char(k,C)$ to itself, then the triple $(C,T,k)$
is {\it splitting rigid}.
\end{prop}

We are thus led to compute the group of homotopy classes of self
homotopy equivalences $Aut(C \times T)$ of $C\times T$.

Before we can do this , let's recall a work by Booth and Heath
\cite{bh}

Given spaces with base point $(X,x_0)$ and $(Y,y_0)$, there is a
natural map $$\varphi:map(X \times Y, X \times Y) \to map(X,X)
\times map(Y,Y)$$ defined by $\varphi(f)=(g,h)$ where $g(x)=\pi_X
\circ f(x,y_0) $ , $h(y)=\pi_Y \circ f(x_0,y) $  and
$\pi_X:X\times X , \pi_Y:X \times Y \to Y$ are projections to the
factors $X$ and $Y$ respectively.
\begin{defn}\label{T:def}
Let $X$ and $Y$ be two spaces . We say $X$ and $Y$ have the
induced equivalence property(IEP) if whenever $f$ is a homotopy
equivalence , then $g,h$ defined above are homotopy equivalences.
\end{defn}
\begin{rem}
Let  $X$ and $Y$ be such that for each $i >0$, at least one of
$\pi_i(X)$ and $\pi_i(Y)$ is zero. Then they satisfy the IEP by
Whitehead theorem.
\end{rem}
With the above notion, we can quote the following
\begin{thm}
Let $X$ and $Y$ be two spaces having IEP. Suppose further that
$[X,map(Y,Y)_{id}]=0$, then there is a short exact sequence of
groups and homomorphisms
\[
1 \to [Y,map(X,X)_{id}] \overset{\theta}{\to} Aut(X \times Y) \to
Aut(X) \times Aut(Y) \to 1
\]
which splits by a homomorphism $\sigma:Aut(X) \times Aut(Y) \to
Aut(X \times Y)$ given by $\sigma(g,h)=g \times h$
\end{thm}

Let $X=C$ and $Y=T$ where $C,T$ be as in Theorem\ref{T:main}. Then
$X$ and $Y$ have IEP by the remark following the
definition\ref{T:def}. On the other hand, $[C,map(T,T)_{id}]=0$
since it is well known that $map(T,T)_{id}=T$ and $C$ is
1-connected and thus first cohomology of $C$ is trivial.

Now given $f \in Aut(C \times T)$ , to prove that the induced
homomorphism in cohomology maps $Char(k,C)$ to itself , it
suffices to assume that $f \in Im([T,map(C,C)_{id}])$ by the exact
sequence above. Recall that the map $\theta:[T,map(C,C)_{id}] \to
Aut(C \times T)$ is given by $\theta(f)(x,y)=(f(y)(x),y)$. The
above argument gives the following
\begin{cor}\label{T:key}
Let $C$ be a closed simply-connected manifold, $T$ be a torus, $k$
be a positive integer. Then the triple $(C,T,k)$ is {\it splitting
rigid} if, for any map $f: T \to map(C,C)_{id}$ , the adjoint
$\tilde{f}: T\times C \to C$  induces a homomorphism in cohomology
given by ${\tilde{f}}^*(u)= 1 \otimes u$ for any $u \in H^*(C) $ .
\end{cor}

\section{The proof of  Theorem\ref{T:main} }

\begin{proof}[Proof of  Theorem\ref{T:main}]
Let $T=(S^1)^s$. By Corollary \ref{T:key}, to prove
Theorem\ref{T:main} , it suffices to prove that , for map $f:
T\times C \to C$ such that $f(y_0,-)$ homotopic to $id$, it
induces a homomorphism in cohomology given by $f^*(u)= 1 \otimes
u$. Given such $f$, for any $u \in H^*(C) $, $$f^*(u)=1 \otimes u
+ \underset{k }{\sum} \underset{i_1  \cdots  i_k
}{\sum}\lambda_{i_1 \cdots i_k}(u) \otimes \iota_{i_1}\cdots
\iota_{i_k}$$ where the first sum is taken over $k$ from $1$ to
$s$ and the second sum is taken over all $(i_1\cdots i_k)'s$ such
that  $0< i_1 < \cdots < i_k < s+1$. Thus we get a sequence of
maps $\lambda_{i_1 \cdots i_k}:H^n(C) \to H^{n-k}(C)$ where $0<k<
s+1$ and $0< i_1 < \cdots < i_k < s+1$.

To prove that $f^*(u)= 1 \otimes u$ , it suffices to prove that
$\lambda_{i_1 \cdots i_k}=0$ where $0<k< s+1$ and $0< i_1 < \cdots
< i_k < s+1$ while for the proof of later we need to study the
behaviour of these maps with respect to the cup product of
cohomology.

If $u ,v \in H^*(C)$ , then $$f^*(uv)=1 \otimes uv + \underset{k
}{\sum} \underset{i_1  \cdots  i_k }{\sum}\lambda_{i_1 \cdots
i_k}(uv) \otimes \iota_{i_1}\cdots \iota_{i_k}$$

 On the other hand $f^*(uv)=f^*(u)f^*(v)$. Using the formula for $f^*(uv),
 f^*(u),f^*(v)$ and comparing the terms associated with $\iota_{i_1}\cdots
\iota_{i_k}$, we find  the following equations
\[
\lambda_{i_1 \cdots i_k}(uv)=\lambda_{i_1 \cdots i_k}(u)v +
(-1)^{k|v|}u\lambda_{i_1 \cdots i_k}(v) \pm \sum \lambda_{j_1
\cdots j_p}(u)\lambda_{l_1 \cdots l_q}(v)\otimes \iota_{i_1}\cdots
\iota_{i_k}
\] where $p+q=k$ with $p>0,q>0$ and the sum is taken over all partitions of
$i_1,\cdots,i_k$ into $j_1<\cdots <j_p$ and $l_1<\cdots <j_q$.

Let $k=1$. Then then above formula implies that $\lambda_{i_1}$ is
a derivation of degree $-1$ which is trivial by the condition of
the Theorem\ref{T:main}. The above formula in case $k=2$ implies
that $\lambda_{i_1i_2}$ is a derivation of degree $-2$ modulo
products of derivations of degree $-1$ which are trivial. It
follows that $\lambda_{i_1i_2}$ is a derivation of degree $-2$ and
thus is trivial by the condition of the Theorem\ref{T:main}.
Inductively we can prove that all $\lambda_{i_1 \cdots i_k}$ are
trivial which completes the proof of the Theorem\ref{T:main}.

\end{proof}
----------------------------------------------------------------

----------------------------------------------------------------

\end{document}